\documentclass[11pt]{amsart}

\usepackage{amsmath}
\usepackage{amssymb}
\usepackage{graphicx}
\usepackage{amscd}
\usepackage{stmaryrd}
\usepackage[mathscr]{euscript}

\theoremstyle{definition}

\numberwithin{equation}{section}

\newcommand{\sB}{\mathscr B}
\newcommand{\sC}{\mathscr C}
\newcommand{\sD}{\mathscr D}
\newcommand{\sE}{\mathscr E}

\newcommand{\sH}{\mathscr H}
\newcommand{\sJ}{\mathscr J}
\newcommand{\sK}{\mathscr K}

\newcommand{\sP}{\mathscr P}
\newcommand{\sR}{\mathscr R}

\newcommand{\sX}{\mathscr X}
\newcommand{\sY}{\mathscr Y}

\newcommand{\p}{\partial}

\title[Miscellaneous on Commutants]{Miscellaneous on Commutants $\mbox{\rm mod}$ Normed Ideals and Quasicentral Modulus $I$}

\author[D.-V. Voiculescu]{Dan-Virgil Voiculescu}

\address{Department of Mathematics \\ University of California at Berkeley \\ Berkeley, CA\ \ 94720-3840}
\email{{\tt dvv@math.berkeley.edu}}

\dedicatory{Dedicated to the memory of Vaughan Jones}

\begin{document}

\begin{center}
\end{center}

\begin{abstract}
We define commutants mod normed ideals associated with compact smooth manifolds with boundary. The results about the $K$-theory of these operator algebras include an exact sequence for the connected sum of manifolds, derived from the Mayer--Vietoris sequence. We also make a few remarks about bicommutants mod normed ideals and about the quasicentral modulus for the quasinormed $p$-Schatten--von~Neumann classes $0 < p < 1$.
\end{abstract}

\maketitle


\section{Introduction}
\label{sec1}

We arrived at the quasicentral modulus, a number associated to a $n$-tuple $\tau$ of operators and a normed ideal $(\sJ,|\quad|_{\sJ})$ of compact operators, studying perturbations of $n$-tuples of operators with perturbations from the normed ideal (\cite{26},\cite{27}). It turned out to play a key role in the invariance of absolutely continuous parts and diagonalizability $\mbox{mod}$ the normed ideal in the case of $n$-tuples of commuting Hermitian operators, and found uses in non-commutative geometry \cite{6}. The number is also of interest in the non-commutative setting of finitely generated groups and it also appears to be related to dynamical entropy. The number is a kind of measure of the $n$-tuple of operators at a scale given by the normed ideal. More recently, trying to understand the ubiquity of the quasicentral modulus, we began looking for more structure around this number and we found that the commutant modulo the normed ideal of the $n$-tuple of operators was relevant structure. Vaguely, this resembles a $K$-theory construction of the commutant modulo the compact operators, used in the Paschke duals, but there are unexpected twists in the analogy, like getting $C^*$-algebras in situations where one only expected Banach algebras or a much richer $K$-theory (see \cite{31} for a recent survey).

\bigskip
\noindent
\underline{\qquad\qquad} 

2020 {\em Mathematics Subject Classification}. Primary: 47A55; Secondary: 46L80, 46L87, 47B10.

{\em Key words and phrases}. commutant mod normed ideal, bicommutant mod normed ideal, $K$-groups of $C^*$-algebras, connected sum of smooth manifolds.


\clearpage
In this paper we have collected various remarks and results further extending the frame we began constructing with the commutants mod normed ideals.

We begin with the shortest discussion of a topic, actually only to point out some questions. Perturbations of Hermitian operators with singular spectral measure, as shown by rank one perturbation facts, behave very differently from what happens in the case of dimension $\ge 1$, where the quasicentral modulus works well. While rank one perturbations certainly don't fit with quasicentral modulus technique one may wonder whether these results signal some more general rough pathologies for perturbation from the quasinormed Schatten--von~Neumann classes with $0 < p < 1$.

Our second topic is completing the picture of the commutant mod a normed ideal, by introducing the bicommutant mod and normed ideal, with the analogy with von~Neumann algebras in mind. Here we can use the results about commutants modulo the compacts of von~Neumann algebras and their refinement to normed ideals \cite{3}, which was preceded by (\cite{14},\cite{15},\cite{23}). What this gives is that the bicommutant mod the normed ideal, is the normed ideal plus a certain ``smooth algebra'' $\sD(\tau,\sJ)$ for the $n$-tuple $\tau$ contained in its von~Neumann algebra. In the commutative setting this can be stated in terms of functional calculus and here the multivariable operator theory work in (\cite{10},\cite{17},\cite{21}) is relevant to getting a grasp on these smoothness classes.

Our last focus, which is the main part of this paper is about commutants mod normed ideals associated to compact differentiable manifolds with boundary. Here a choice of an embedding into some ${\mathbb R}^m$ provides an $m$-tuple of smooth functions the multiplication operators by which determine the commutant mod the normed ideal up to equivalent norms. We use this to develop for these commutants mod several of the results we had about the case of $n$-tuples of operators. We also show that the $K$-theory of Calkin algebras for commutants mod normed ideals for smooth manifolds has a Mayer Vietoris sequence when we perform connected sums on the manifolds. The key to such result is that this kind of coronas/``Calkin algebras'' are $C(\sX)$ algebras in the sense of Kasparov or more general $\sX$-algebras considered by Kirchberg. The conclusion from these results seems to be that operations like connected sums may be manageable, but the main problem to better understand these commutants mod for manifolds would be with understanding better their structure or $K$-theory in simple cases like balls or spheres.

The paper has six sections, the first being this introduction and the last the list of references. The second section is about preliminaries. Then the third section is about the problems concerning the quasinormed setting, the fourth about bicommutants mod the fifth about commutants mod for differentiable manifolds with boundary.

\section{Preliminaries}
\label{sec2}

By $\sH$ we will denote a separable infinite-dimensional complex Hilbert space and $\sB(\sH)$, $\sK(\sH)$, $\sB/\sK(\sH)$, $\sR(\sH)$, $\sP(\sH)$ or simply $\sB$, $\sK$, $\sB/\sK$, $\sR$, $\sP$ will denote the bounded operators, the compact operators, the Calkin algebra, the finite rank operators and the finite rank Hermitian projectors. The canonical homomorphism $\sB \to \sB/\sK$ will be denoted by $p$.

The normed ideals of compact operators $(\sJ,|\quad|_{\sJ})$ or $(\sJ(\sH),|\quad|_{\sJ})$, when $\sH$ needs to be specified, and their properties are discussed in the standard reference (\cite{11},\cite{24}). Note that the normed ideals, as defined in these references have the so-called Fatou property which is required in results of \cite{3}, which we shall use. In particular $(\sC_p,|\quad|_p)$ will denote the Schatten--von~Neumann $p$-class, $1 \le p < \infty$. We shall also use the notation $(\sC_p,|\quad|_p)$ when $0 < p < 1$, in which case $|\quad|_p$ is a quasi-norm and the ideal $\sC_p$ is no longer locally convex. When $1 \le p \le \infty$ we denote by $(\sC^-_p,|\quad|_p^-)$ the $(p,1)$ Lorentz ideal, the norm on which is given by 
\[
|T|_p^- = \sum_{j \in {\mathbb N}} s_jj^{-1+1/p}
\]
where $s_1 \ge s_2 \ge \dots$ are the eigenvalues of $(T^*T)^{1/2}$. Note that $\sC_1^- = \sC_1$, $|\quad|_1^- = |\quad|_1$.

If $\tau = (T_j)_{1 \le j \le n}$ is a $n$-tuple of bounded operators, we denote by $\tau^* = (T_j^*)_{1 \le j \le n}$ $\|\tau\| = \max_{1 \le j \le n} \|T_j\|$ and by $[\tau,X] = ([T_j,X])_{1 \le j \le n}$  where $X \in \sB$. If $\sigma = (S_k)_{1 \le k \le m}$ is an $m$-tuple of bounded operators then $(\sigma,\tau)$ is the $(m+n)$-tuple $(S_1,\dots,S_m,T_1,\dots,T_n)$ and if $m = n$, $\sigma + \tau = (S_j + T_j)_{1 \le j \le n}$. If $(\sJ,|\quad|_{\sJ})$ is a normed ideal,
\[
\sE(\tau;\sJ) = \{X \in \sB \mid [T_j,X] \in \sJ, 1 \le j \le n\}
\]
is the commutant of $\tau \mod \sJ$. $\sE(\tau,\sJ)$ is a Banach algebra with the norm
\[
\||X\|| = \|X\| + |[\tau,X]|_{\sJ}
\]
where the $|\quad|_{\sJ}$-norm is the max of the $\sJ$-norms of the components. If $\tau = \tau^*$, $\sE(\tau;\sJ)$ is a $*$-algebra and $\||X\|| = \||X^*\||$. We also denote by
\[
\begin{aligned}
\sK(\tau;\sJ) &= \sE(\tau,\sJ) \cap \sK \\
\sE/\sK(\tau,\sJ) &= \sE(\tau;\sJ)/\sK(\tau;\sJ) \\
\pi: \sE(\tau;\sJ) &\rightarrow \sE/\sK(\tau;\sJ)
\end{aligned}
\]
the compact ideal, the quotient Banach algebra and the canonical homomorphism. Note that $p(\sE(\tau;\sJ)) \in \sB/\sK$ and $\sE/\sK(\tau;\sJ)$ are algebraically isomorphic and the map
\[
\sE/\sK(\tau;\sJ) \to p(\sE(\tau;\sJ))
\]
is contractive. The quasicentral modulus of $\tau$ w.r.t. $\sJ$ is the number
\[
k_{\sJ}(\tau) = \liminf_{A \in \sR^+_1} |[A,\tau]|_{\sJ}
\]
where $\sR^+_1 = \{X \in \sR \mid 0 \le X \le I\}$ is the set of finite rank positive contractions endowed with the natural order (\cite{26}, \cite{28}). If $\sJ = \sC_p$ or $\sJ = \sC^-_p$ we denote $k_{\sJ}(\tau)$ also by $k_p(\tau)$ and respectively $k^-_p(\tau)$. In case $0 < p < 1$ the definition of $k_p(\tau)$ still makes sense with $|\quad|_p$ being the $p$-quasinorm in this case. We shall also consider the modulus of quasidiagonality
\[
qd_{\sJ}(\tau) = \lim_{P \in \sP} \inf|[P,\tau]|_{\sJ}
\]
the liminf being with respect to the natural order on the finite rank projections $\sP$.

Assuming $\sR$ is dense in $\sJ$, we recall (see \cite{23}) that if $k_{\sJ}(\tau) = 0$ then $\sE/\sK(\tau;\sJ)$ is a $C^*$-algebra, while if only $k_{\sJ}(\tau) < \infty$ then $\sE/\sK(\tau;\sJ)$ is isomorphic as a Banach algebra with involution to a $C^*$-algebra (not isometrically), actually $p(\sE(\tau;\sJ))$ is a $C^*$-subalgebra of $\sB/\sK$ and is isomorphic (not isometrically) to $\sE/\sK(\tau;\sJ)$.

We shall also denote by $\sP\sE(\tau;\sJ)$ the Hermitian projectors in $\sE(\tau;\sJ)$.

Finally if $\omega \subset \sB$ we denote by $\omega'$ its commutant in $\sB$
\[
\omega' = \{X \in \sB \mid [X,\omega] = \{0\}\}.
\]
Then $\omega'' = (\omega')'$ will be the bi-commutant of $\omega$, that is in case $\omega = \omega^*$, the von~Neumann algebra generated by $\omega$.

\section{The problem with $0 < p < 1$}
\label{sec3}

A $n$-tuple $\tau$ of commuting Hermitian operators is diagonalizable $\mod \sC_n^-$, which is also equivalent to $k^-_n(\tau) = 0$, iff the spectral measure of $\tau$ is singular w.r.t. $n$-dimensional Lebesgue measure \cite{26}. A similar result also holds for the normed ideal $\sC^-_p$, $1 < p < n$ when the spectrum of $\tau$ is contained in certain fractals of Hausdorff dimension $p$, with $p$-dimensional Hausdorff measure replacing Lebesgue measure (\cite{8},\cite{32}). In these cases we observe that diagonalizability $\mod \sC^-_p$ and $k^-_p$ are tied to $p$-Hausdorff measure, where $p \ge 1$.

$p$-Hausdorff measure and dimension are also natural and interesting when $0 < p < 1$. Something quite different occurs with the diagonalization problem and we will not even need the refined Lorentz double scale to see it. If $T = T^*$ has a cyclic vector and its spectral measure is singular w.r.t. Lebesgue measure, then (\cite{9}, \cite{12}, \cite{22}, \cite{25}) there is a rank one operator $X = X^*$ of arbitrary small norm so that $X + T$ is diagonal. A general Hermitian operator $T$ on a separable Hilbert space is an orthogonal sum of a sequence of Hermitian operators with cyclic vector. We can find for each summand a rank one perturbation which diagonalizes it and so that the orthogonal sum of the perturbations is a Hermitian operator $X$ so that $X \in \sC_p$ for all $0 < p < 1$. This means that $k_1(T) = 0$, which is equivalent to the singularity of the spectral measure of $T$ implies $k_p(T) = 0$ for all $0 < p < 1$. Here $k_p(T)$ is defined the same way as in the case $p \ge 1$, with the only difference that the $p$-norm is now a $p$-quasinorm when $0 < p < 1$. This looks like a kind of ``phase-transition'' at $p = 1$ in the correspondence between $k^-_p$ and $p$-Hausdorff measure. It is natural to ask how general this is.

\medskip
\noindent
{\bf Problem 3.1.} Let $\tau$ be a $n$-tuple of Hermitian operators. Assume $k_1(\tau) = 0$. Does it follow that $k_p(\tau) = 0$ for $0 < p < 1$? If this is too general, one can add the assumption that $\tau$ is a $n$-tuple of commuting Hermitian operators.

\medskip
\noindent
{\bf Remark 3.1.} If $T = T^*$ then $k_1(T) = 0$ is equivalent to $qd_1(T) = 0$ implies $k^-_p(\tau) = 0$ or even $qd_p(\tau) = 0$ for $0 < p < 1$. Again one may consider the extra assumption about commuting components in $\tau$.

\section{The bicommutant $\mbox{\rm mod}$ a normed ideal}
\label{sec4}

In this section $\tau = (T_j)_{1 \le j \le n}$ is a $n$-tuple of Hermitian operators on $\sH$ and $(\sJ,|\quad|_{\sJ})$ is a normed ideal.

\medskip
\noindent
{\bf Definition 4.1.} The bicommutant of $\tau \mod \sJ$ is the $*$-algebra
\[
\sE\sE(\tau;\sJ) = \{X \in \sB \mid [X,\sE(\tau;\sJ)] \subset \sJ\}.
\]
Its compact ideal is
\[
\sK\sE(\tau;\sJ) = \sE\sE(\tau;\sJ) \cap \sK
\]
and the bicommutant Calkin algebra is
\[
\sE/\sK\sE(\tau;\sJ) = \sE\sE(\tau;\sJ)/\sK\sE(\tau;\sJ).
\]

\medskip
\noindent
{\bf Lemma 4.1.} {\em If $X \in \sE\sE(\tau;\sJ)$ then
\[
\sup\{|[X,Y]|_{\sJ} \mid Y \in \sE(\tau;\sJ),\||Y\|| \le 1\} < \infty
\]
and $\sE\sE(\tau;\sJ)$ is an involutive Banach algebra with isometric involution when endowed with the norm
\[
\|\|X\|\| = \|X\| + \sup\{|[X,Y]|_{\sJ} \mid Y \in \sE(\tau;\sJ),\||Y\|| \le 1\}.
\]
}

\medskip
The lemma is a consequence of the closed graph theorem applied to the map
\[
\sE(\tau;\sJ) \ni Y \to [X,Y] \in \sJ
\]
where $X \in \sE\sE(\tau;\sJ)$. Indeed if $\||Y_n-Y\|| \to 0$, $|[X,Y_n]-Z|_{\sJ} \to 0$ then $[X,Y] = Z$, so the  graph is closed. The rest is also an easy exercise.

\medskip
\noindent
{\bf Theorem 4.1.} {\em We have $\sE\sE(\tau;\sJ) = \sE\sE(\tau;\sJ) \cap (\tau)'' + \sJ$ and $\sE\sE(\tau;\sJ) \cap (\tau)''$ is a closed subalgebra of $\sE\sE(\tau;\sJ)$.
}

\medskip
\noindent
{\bf Proof.} We have $(\tau)' \subset \sE(\tau;\sJ)$ so that if $X \in \sE\sE(\tau;\sJ)$ then we will have $[X,(\tau)'] \subset \sJ$. By \cite{3} this implies $X \in (\tau)'' + \sJ$. On the other hand $\sJ \subset \sE\sE(\tau;\sJ)$ so this gives
\[
\sE\sE(\tau;\sJ) = \sE\sE(\tau;\sJ) \cap (\tau)'' + \sJ.
\]
If $Y \in \sE(\tau;\sJ)$, the map
\[
\sE\sE(\tau;\sJ) \ni X \to [X,Y] \in \sJ
\]
is continuous. Since $(\tau)' \subset \sE(\tau;\sJ)$ it follows that $\sE\sE(\tau;\sJ) \cap (\tau)''$ is closed in $\sE\sE(\tau;\sJ)$.\qed

\medskip
\noindent
{\bf Corollary 4.1.} {\em 
We have $p(\sE\sE(\tau;\sJ)) \subset p((\tau)'') \subset \sB/\sK$. If $(\tau)'' \cap \sK = \{0\}$, then $p$ gives an isomorphism of $\sE\sE(\tau;\sJ) \cap (\tau)''$ and $p(\sE\sE(\tau;\sJ))$.
}

\medskip
The corollary is an immediate consequence of the theorem.

\medskip
\noindent
{\bf Definition 4.2.} We define $\sD(\tau;\sJ) = \sE\sE(\tau;\sJ) \cap (\tau)''$ to be the {\em smooth algebra of $\tau$ w.r.t. $\sJ$}. We also define $\sD/\sK(\tau;\sJ) = p(\sD(\tau;\sJ)) = p(\sE\sE(\tau;\sJ))$ to be the {\em essential smooth algebra of $\tau$ w.r.t. $\sJ$}.

\medskip
\noindent
{\bf Remark 4.1.} If $\tau-\tau' \in (\sJ)^n$ then clearly $\sE(\tau;\sJ) = \sE(\tau';\sJ)$ so that $\sE\sE(\tau;\sJ) = \sE\sE(\tau';\sJ)$. While $(\tau)''$ and $(\tau')''$ are different in general we have obviously $\sD/\sK(\tau;\sJ) = \sD/\sK(\tau';\sJ)$. Actually even more $\sE\sE(\tau';\sJ)/\sJ = \sE\sE(\tau;\sJ)/\sJ$ gives that $\sD(\tau;\sJ)$ and $\sD(\tau';\sJ)$ have the same image in $\sB/\sJ$.

\medskip
\noindent
{\bf Remark 4.2.} Note that if $\tau$ is an $n$-tuple of commuting Hermitian operators, then $(\tau)''$ is the algebra of Borel functions of $\tau$ by functional calculus and this relates $\sD(\tau;\sJ)$ to smoothness properties of operator functions, see for instance \cite{10}, \cite{17}, \cite{21} and related work.

\section{Commutants ${\rm mod}$ normed ideals associated with differentiable manifolds}
\label{sec5}

Let $\sX$ be a compact $C^\infty$-manifold with boundary (\cite{13}). Let $\mu$ be a Radon measure on $\sX$ such that in each local parametrization the restriction of $\mu$ has the same absolute continuity class as Lebesgue measure. The representation $M: L^{\infty}(\sX;\mu) \to \sB(L^2(\sX;\mu)$ by multiplication operators $M(f)h = fh$, up to unitary equivalence does not depend on the choice of $\mu$. The intertwining operator for the representations arising from two such measures is given by multiplication with the square-root of their Radon--Nikodym derivative. Let $(\sJ,|\ |_{\sJ})$ be a normed ideal, which to simplify, is so that $\sR$ is dense in $\sJ$. Then we define
\[
\sE(\sX,\sJ) = \{X \in \sB(L^2(\sX,\mu)) \mid [M(C^{\infty}(\sX)),X] \subset \sJ\}
\]
the $*$-algebra which is the commutant $\mod \sJ$ of $M(C^{\infty}(\sX))$. If $\alpha$ is a diffeomorphism of $\sX$ there is a unitary operator $U_{\alpha}h = (h \circ \alpha^{-1}) \cdot [\alpha,\mu : \mu]^{1/2}$ on $L^2(\sX,\mu)$ so that $\alpha \to U_{\alpha}$ is a representation of $\sD${\em iff} $(\sX)$ and
\[
U_{\alpha}M(f)U_{\alpha}^{-1} = M(f \circ\alpha^{-1}).
\]
It follows that
\[
U_{\alpha}\sE(\sX;\sJ)U_{\alpha}^{-1} = \sE(\sX;\sJ)
\]
and we see that
\[
\alpha \to (X \to U_{\alpha}XU_{\alpha}^{-1})
\]
gives a homomorphism of $\sD${\em iff} $(\sX)$ into the automorphisms of $\sE(\sX;\sJ)$. Our next aim is to show that $\sE(\sX;\sJ)$ can actually be defined using only finitely many $M(f)$ so that we can use the results on commutants mod a normed ideal of commuting $n$-tuples of Hermitian operators. In particular $\sE(\sX;\sJ)$ is a Banach algebra with involution.

\medskip
\noindent
{\bf Lemma 5.1.} {\em Let $\tau = (X_j)_{1\le j \le n}$ be a $n$-tuple of commuting Hermitian operators and $\Omega \subset {\mathbb R}^n$ an open set so that the spectrum $\sigma(\tau) \subset \Omega$. If $\sJ$ is a normed ideal and $F: \Omega \to {\mathbb R}$ is a $C^{\infty}$-function then if $Y \in \sB(L^2(\sX , \mu))$ then
\[
|[F(\tau),Y]|_{\sJ} \le C|[\tau,Y]|_{\sJ}
\]
where the constant does not depend on $\sJ$.}

\medskip
This is well-known. It can for instance be proved by dealing first with the case of $F(x_1,\dots,x_n) = \exp(i \sum_{1 \le j \le n}x_j\xi_j)$ and then obtaining the result for general $F$ by using the fact that the Fourier transform can be used to express $F(\tau)$ as an integral involving exponentials applied to $\tau$.

\medskip
\noindent
{\bf Proposition 5.1.} {\em Let $f_1,\dots,f_m \in C^{\infty}(\sX)$ be such that $F: \sX \to {\mathbb R}^m$, where $F = (f_1,\dots,f_m)$ is a neat embedding of $\sX$ into a half-space of ${\mathbb R}^m$, in the sense of \cite{13}. Then we have
\[
\sE(\sX;\sJ) = \sE((M(f_j))_{1\le j \le m};\sJ).
\]
}

\begin{proof}
Obviously we have $LHS \subset RHS$ so what we need to prove is that $LHS \supset RHS$. We need to prove that if $X \in \sE((M(f_j))_{1 \le j \le m};\sJ)$ and $g \in C^{\infty}(\sX)$ is real-valued then $[M(g),X] \in \sJ$. This follows from Lemma~$5.1$ since we can find a $C^{\infty}$-function $G: \Omega \to {\mathbb R}$, where $\Omega$ is a neighborhood of $F(\sX)$, so that $G \circ F = g$ and then $G(M(f_1),\dots,M(f_m)) = M(g)$.
\end{proof}

Since $\sX$ is compact, by (\cite{13}, Thm.~$4.3$) the embedding assumption in the preceding proposition can be satisfied and thus the results we have for commutants mod normed ideals of $m$-tuples of Hermitian operators apply to $\sE(\sX;\sJ)$. In what follows we list a few of the consequences which we would like to point out in particular. An embedding of $\sX$ provides a norm on $\sE(\sX;\sJ)$ w.r.t.\ which it is a Banach algebra with isometric involution and using the closed graph theorem one finds that different embeddings give rise to equivalent norms on $\sE(\sX;\sJ)$. Then $\sE(\sX;\sJ) \cap \sK$ which we shall denote by $\sK(\sX;\sJ)$ is a closed ideal and $\sE/\sK(\sX;\sJ) = \sE(\sX;\sJ)/\sK(\sX;\sJ)$ is also a Banach algebra with involution.

\medskip
\noindent
{\bf Theorem 5.1.} {\em 
Let $\sX$ be a compact smooth manifold with boundary and real-valued $f_1,\dots,f_m \in C^{\infty}(\sX)$ which give a neat imbedding of $\sX$ (\cite{13}). Let further $(\sJ,|\ |_{\sJ})$ be a normed ideal so that $\sJ \supset \sC^-_n$, $\sJ \ne \sC^-_n$ and with $\sR$ dense in $\sJ$ where $n$ is the dimension of $\sX$. 

Then:

{\rm a)} $k_{\sJ}(M(f_1),\dots,M(f_m)) = 0$, $0 < k^-_n(M(f_1),\dots,M(f_m)) < \infty$ and if $\sH_1 \subset L^2(M,\mu)$ is an invariant subspace for $(M(f_1),\dots,M(f_m))$, $\sH_1 \ne 0$, then still $k^-_n(M(f_1)|_{\sH_1},\dots,M(f_m)|_{\sH_1}) > 0$.

{\rm b)} $\sE/\sK(\sX;\sC^-_n)$ is (non-isometrically) isomorphic to a $C^*$-algebra and $\sE/\sK(\sX;\sJ)$ is a $C^*$-algebra, the norm being independent of the embedding.

{\rm c)} The centre of $\sE/\sK(\sX;\sC^-_n)$ is $\pi(M(C(X)))$. Also the centre of $\sE/\sK(\sX;\sJ)$ is $\pi(M(C(\sX)))$.

{\rm d)} If $n \ge 3$ and $v_m \in C^{\infty}(\sX)$ are so that $|v_m| = 1$ and $w - \lim_{m \to \infty} M(v_m) = 0$ then the strong limit
\[
\Phi(T) = {\underset{m \to \infty}{s-\lim}} M(v_m)TM(v_m)^*
\]
exists for every $T \in \sE(\sX;\sC^-_n)$ and is independent of the choice of $v_m$'s and gives a $*$-homomorphism $\Phi:\sE(\sX;\sC^-_n) \to M(L^{\infty}(\sX,\mu))$ so that $\Phi(M(f)) = M(f)$, if $f \in 
L^{\infty}(\sX,\mu)$, $\|\Phi\| = 1$ and $\Phi(\sK(\sX,\sC^-_n)) = \{0\}$. In particular there is a $*$-homomorphism $\Psi: \sE/\sK(\sX;\sC^-_n) \to \pi(M(L^{\infty}(\sX;\mu))$ so that $\pi \circ \Phi = \Psi \circ \pi$.
}

\begin{proof}
a) One can cover $\sX$ with a finite number of open sets in each of which there is some subset of 
$\{f_1,\dots,f_m\}$ which provides a coordinate system. Then $L^2(\sX;\mu)$ is an orthogonal sum of $L^2(\Omega;\mu\mid\Omega)$ where $\Omega$ is a Borel subset of one of those open sets. Then the assertion reduces to that for an $n+p$-tuple of commuting Hermitian operators the first $n$ of which have joint spectrum which is Lebesgue of multiplicity one and the remaining $p$ are smooth functions of these. Thus a) becomes a consequence of the facts about normed ideal perturbations of commuting $n$-tuples of Hermitian operators (\cite{4}, \cite{26}, \cite{27}, \cite{28}) supplemented by Lemma~$5.1$.

b) This is a consequence of a) and of \cite{29}.

c) This follows from \cite{30} and a).

d) Let $v_m \in C^{\infty}(\sX)$ be so that $|v_m| = 1$ and $w - \lim_{m \to \infty} M(v_m) = 0$, then the strong limit
\[
s - \lim_{m \to \infty} M(v_m)TM(v_m)^* = \Phi(T)
\]
exists for every $T \in \sE(\sX;\sC^-_n)$. This is obtained after a construction like in the proof of a) from \cite{27} Thm.~$2.3$, the comments following it and Corollary~$1.6$ in \cite{27}. As remarked in \cite{27} the fact that the sequence of $v_m$ is arbitrary implies, $\Phi$ does not depend on the choice of $v_m$'s which implies $\Phi(T) \in (M(C^{\infty}(\sX)))' = M(L^{\infty}(\sX;\mu))$. Note that if $T \in \sK(\sX,\sC^-_n)$ then $\Phi(T) \in \sK$ and $\sK \cap M(L^{\infty}(\sX;\mu)) = \{0\}$ gives $\Phi(T) = 0$.
\end{proof}

It is also useful to make the following simple observation. If we compress $\sE(\sX;\sJ)$ to $L^2(\Omega;\mu)$ where $\Omega \subset \sX$ is a Borel set, we get the commutant $\mod \sJ$ of the multiplication operators by $f|_{\Omega}$ where $f \in C^{\infty}(\sX)$ or equivalently of $f_1|_{\Omega},\dots,f_n|_{\Omega}$ for an embedding $f_1,\dots,f_m$ of $\sX$. Note also these kind of considerations carry over to $\sK(\sX;\sJ)$ and $\sE/\sK(\sX;\sJ)$. In particular we have the following, for an arbitrary normed ideal $\sJ$.

\medskip
\noindent
{\bf Lemma 5.2.} {\em
Let $A \subset \sX\backslash\partial\sX$ be a closed submanifold with boundary of the same dimension as $\sX$ and $P = M(\sX_A)$. Then
\[
\sE(A;\sJ) = P\sE(\sX;\sJ) \mid L^2(A,\mu\mid A)
\]
and we also have isometric isomorphisms of $\sE/\sK(A;\sJ)$ and $\pi(P)\sE/\sK(\sX;\sJ)\pi(P)$.
}

\medskip
If $\sJ \supseteq \sC^-_n$ by c) of Theorem~5.1 the center of $\sE/\sK(\sX;\sJ)$ being $\pi(M(C(\sX)))$ to an open set $U \subset \sX$ there corresponds a closed two-sided ideal $\sE/\sK(U,\sX;\sJ)$ of $\sE/\sK(\sX;\sJ)$ which can be defined as follows. If $A \subset \sX$ is a Borel set let $P_A$ be the projection $M(\sX_A)$ where $\sX_A$ is the indicator function. We define $\sE/\sK(U,\sX;\sJ)$ to be the closure of
\[
\bigcup_{K \subset\subset U} \pi(P_K)\sE/\sK(\sX;\sJ)\pi(P_K)
\]
where $K \subset\subset U$ denotes the fact that $K \subset U$ is a compact subset.

\medskip
\noindent
{\bf Lemma 5.3.} {\em
{\rm a)} If $K_j \subset\subset U$, $j \in {\mathbb N}$ are so that $K_j \subset \overset{\circ}{K}_{j+1}$ and $\bigcup_{j \in {\mathbb N}} K_j = U$ then $\sE/\sK(U,\sX;\sJ)$ is the closure of $\bigcup_{j \in {\mathbb N}} \pi(P_{K_j})\sE/\sK(\sX;\sJ)\pi(P_{K_j})$.

{\rm b)} If $f \in C(\sX)$ is so that $f^{-1}(0) = \sX\backslash U$ then $\sE/\sK(U,\sX;\sJ)$ is the closure of
\[
\pi(M(f))\sE/\sK(\sX;\sJ).
\]

{\rm c)} $\sE/\sK(U,\sX;\sJ)$ is the closed two-sided ideal of $\sE/\sK(\sX;\sJ)$ generated by $\{\pi(M(f)) \mid f \in C(\sX),f \mid \sX\backslash U = 0\}$.
}

\begin{proof}
a) is obvious since every compact set $K \subset\subset U$ is contained in some $K_j$.

b) Let $K_j = \{x \in \sX \mid |f(x)| \ge 1/j\}$. Since $K_j \subset\subset U$ and $\|| P_{K_j}M(f) - M(f)\|| \to 0$ we get easily $\pi(M(f))\sE/\sK(\sX;\sJ) \subset \sE/\sK(U,\sX;\sJ)$. Also \linebreak
$\pi(M(f))\sE/\sK(\sX;\sJ) \supset \pi(M(f))\pi(P_{K_j})\sE/\sK(\sX;\sJ)\pi(P_{K_j}) \supset \pi(P_{K_j})\sE/\sK(\sX;\sJ)\pi(P_{K_j})$ and a) gives the converse inclusion.

c) Since $\pi(M(f))$ is in the center of $\sE/\sK(\sX;\sJ)$ we have that $\pi(M(f))\sE/\sK(\sX;\sJ)$ is a two-sided ideal and we can also use b) to get c).
\end{proof}

We shall look at the case when the open set $U = \sX\backslash N$ where $N$ is a submanifold satisfying certain conditions so that we can use tubular neighborhood results (\cite{13}) when describing $\sE/\sK(U,\sX;\sJ)$. We shall assume $N$ is a compact submanifold without boundary in $\sX$ (thus $N \cap \partial X = \emptyset$). Using a tubular neighborhood of $N$ and an orthogonal structure on the vector bundle from which it arises (see \cite{13}, p.~116--117) the open disk-subbundle of radius $\epsilon > 0$ yields an open neighborhood $N_{\epsilon}$ of $N$ so that the closure ${\bar N}_{\epsilon}$ is a compact submanifold of $\sX$ with boundary $\p N_{\epsilon} \subset \sX\backslash\partial\sX$. Note that $K_{\epsilon} = \sX \backslash N_{\epsilon}$ is then a compact manifold with boundary $\partial N_{\epsilon} \cup \partial \sX$. Since $0 < \epsilon_1 < \epsilon_2 \Rightarrow N_{\epsilon_1} \subset N_{\epsilon_2}$ and $\bigcap_{\epsilon > 0} N_{\epsilon} = N$ we can use Lemma~$5.3$ a) which gives that if $\sJ \supseteq \sC^-_n$ $\sE/\sK(\sX\backslash N,\sX;\sJ)$ is the closure of $\bigcup_{\epsilon > 0} \pi(P_{\epsilon})\sE/\sK(\sX;\sJ)\pi(P_{\epsilon})$ where $P_{\epsilon} = M(\chi_{\xi\backslash N_{\epsilon}})$. Of course, the receding union can be taken over $\epsilon \in \{1/j \mid j \in {\mathbb N}\}$. Then by Lemma~$5.2$, $\pi(P_{\epsilon})\sE/\sK(\sX;\sJ)\pi(P_{\epsilon})$ identifies with $\sE/\sK(\sX\backslash N_{\epsilon};\sJ)$. Remark also that the assumption $\sJ \supseteq \sC^-_n$ guarantees that these are $C^*$-algebras or at least closed subalgebras of a Banach algebra isomorphic to a $C^*$-algebras. This gives that $\sE/\sK(\sX\backslash N,\sX;\sJ)$ is isomorphic to the inductive limit of the $C^*$-algebras which are isomorphic to the $\sE/\sK(\sX\backslash N_{1/j};\sJ)$. We record this as the next lemma.

\medskip
\noindent
{\bf Lemma 5.4.} {\em Let $N$ be a compact submanifold without boundary of $\sX$ and assume $\sJ \supseteq \sC^-_n$. Then the ideal $\sE/\sK(\sX\backslash N,\sX;\sJ)$ is a Banach algebra which is isomorphic to the inductive limit of the $C^*$-algebras isomorphic to $\sE/\sK(\sX\backslash N_{1/j};\sJ)$ for the $1/j$-neighborhoods of $N$ constructed above.}

\medskip
Next we shall take a look at the effect of connected sums of smooth manifolds on the associated algebras $\sE/\sK(\sX;\sJ)$ when $\sJ \supseteq \sC^-_n$.

Let $\sX,\sY$ be two compact smooth manifolds of the same dimension $n$ and $\sX\#\sY$ their connected sum. Then the construction of the connected sum gives two points $x \in \sX\backslash \partial\sX$, $y \in \sY\backslash \partial\sY$, open sets $U,V \subset \sX\#\sY$ so that $U \cup V = \sX\#\sY$ and diffeomorphisms $\alpha: U \to \sX\backslash\{x\}$, $\beta: V \to \sY\backslash\{y\}$, $\gamma: U \cap V \to (0,1) \times S^{n-1}$. Moreover there are continuous functions $f,g: \sX\#\sY \to [0,\infty)$ so that $f^{-1}((0,\infty)) = U$, $f^{-1}((0,\infty)) = V$, $f + g = 1$ and $pr_1 \circ \gamma = f \mid U \cap V$.

\medskip
\noindent
{\bf Lemma 5.5.} {\em We have
\[
\sE/\sK(U,\sX\#\sY;\sJ) + \sE/\sK(V,\sX\#\sY;\sJ) = \sE/\sK(\sX\#\sY;\sJ)
\]
and
\[
\sE/\sK(U,\sX\#\sY;\sJ) \cap \sE/\sK(V,\sX\#\sY;\sJ) = \sE/\sK(U\cap V,\sX\#\sY;\sJ).
\]
}

\begin{proof}
Lemma~$5.3.$ b) gives
\[
\pi(M(f))\sE/\sK(\sX\#\sY;\sJ) \subset \sE/\sK(U,\sX\#\sY;\sJ)
\]
\[
\pi(M(g))\sE/\sK(\sX\#\sY;\sJ) \subset \sE/\sK(V,\sX\#\sY;\sJ)
\]
and we have $f+g=1$ which yields the first equality. Since $\|f^{1/k}f-f\| \to 0$ and $\|g^{1/k}g-g\| \to 0$ as $k \to \infty$ we see that Lemma~$5.3$ b) implies that $(\pi(M(f^{1/k}))_{k \ge 1}$ and $(\pi(M(g^{1/k})))_{k\ge 1}$ are approximate units of the ideals $\sE/\sK(U,\sX\#\sY;\sJ)$ and $\sE/\sK(V,\sX\#\sY;\sJ)$ respectively. Then $(\pi(M(f^{1/k}))\pi(M(g^{1/k})))_{k\ge 1}$ is an approximate unit for $\sE/\sK(U,\sX\#\sY;\sJ) \cap \sE/\sK(V,\sX\#\sY;\sJ)$. Again by Lemma~$5.3$ b) we have $\pi(M(f^{1/k}))\pi(M(g^{1/k})) = \pi(M((fg)^{1/k}) \in \sE/\sK(U \cap V,\sX\#\sY;\sJ)$ so that $\sE/\sK < (U,\sX\#\sY;\sJ) \cap \sE/\sK(V,\sX\#\sY;\sJ) \subset \sE/\sK(U \cap V,\sX\#\sY;\sJ)$. The opposite inclusion is obvious.
\end{proof}

With these preparations we apply the Mayer--Vietoris exact sequence \cite{12} to this situation.
\[
\scriptsize{\begin{array}{c}
K_0(\sE/\sK(U \cap V,\sX\#\sY;\sJ)) \to K_0(\sE/\sK(U,\sX\#\sY;\sJ)) \oplus K_0(\sE/\sK(V,\sX\#\sY;\sJ)) \to K_0(\sE/\sK(\sX\#\sY;\sJ)) \\
\uparrow\hskip 3.25 in \downarrow \\
K_1(\sE/\sK(\sX\#\sY;\sJ)) \rightarrow K_1(\sE/\sK(U,\sX\#\sY;\sJ)) \oplus K_1(\sE/\sK(V,\sX\#\sY;\sJ)) \rightarrow K_1(\sE/\sK(U \cap V,\sX\#\sY;\sJ))
\end{array}}
\]
We will work more on identifying the terms of this exact sequence so that it relates the groups $K_j(\sE/\sK(\sX\#\sY;\sJ))$ to groups which no longer involve $\sX\#\sY$.

\medskip
\noindent
{\bf Lemma 5.6.} {\em With the notation introduced and the assumption $\sJ \supseteq \sC^-_n$, the diffeomorphisms $\alpha,\beta,\gamma$ give rise to isomorphisms
\[
\sE/\sK(U,\sX\#\sY,\sJ) \simeq \sE/\sK(\sX\backslash \{x\},\sX;\sJ)
\]
\[
\sE/\sK(V,\sX\#\sY,\sJ) \simeq \sE/\sK(\sY\backslash \{y\},\sY;\sJ)
\]
\[
\sE/\sK(U \cap V,\sX\#\sY;\sJ) \simeq \sE/\sK((0,1) \times S^{n-1},[0,1] \times S^{n-1};\sJ)
\]
}

\begin{proof}
We can apply Lemma~$5.4$ to $\sE/\sK(\sX\backslash \{x\},\sX;\sJ)$ with $N = \{x\}$ and the tubular neighborhood arising from the local coordinates chosen for the connected sum gluing. Let $K_j = \sX\backslash N_{1/j}$ and $\alpha^{-1}(K_j) = L_j \subset U$. Then $\sE/\sK(\sX\backslash \{x\},\sX;\sJ)$ is the inductive limit of $\sE/\sK(K_j;\sJ)$ and $\sE/\sK(U,\sX\#\sY;\sJ)$ is the inductive limit of the $\sE/\sK(L_j;\sJ)$. The diffeomorphism $\alpha$ produces an isomorphism of the inductive limits. The identification of $\sE/\sK(V,\sX\#\sY;\sJ)$ and $\sE/\sK(\sY\backslash \{y\},\sY;\sJ)$ using $\beta$ is clearly completely analogous.

For the last assertion let now $K_j = [1/2j,1-1/2j] \times S^{n-1}$ and $L_j = \gamma^{-1}(K_j)$. Then $K_j$ and $L_j$ are compact submanifolds with boundary of $[0,1] \times S^{n-1}$ and $\sX\#\sY$, respectively, which don't intersect the boundaries of these. Using Lemma~$5.2$ and Lemma~$5.3$ we get that $\sE/\sK(U \cap V,\sX\#\sY;\sJ)$ and $\sE/\sK(U \cap V,\sX\#\sY;\sJ)$ are then identified with the inductive of the $\sE/\sK(K_j;\sJ)$ and respectively of the $\sE/\sK(L_j;\sJ)$ and these are isomorphic via isomorphisms induced by $\gamma$.
\end{proof}

With the identifications provided by the preceding Lemma, the reader can further pursue this and work out what the inclusions 
$\sE/\sK(U \cap V,\sX\#\sY;\sJ) \subset \sE/\sK(U,\sX\#\sY;\sJ)$, 
$\sE/\sK(U \cap V,\sX\#\sY;\sJ) \subset \sE/\sK(V,\sX\#\sY;\sJ)$, 
$\sE/\sK(U,\sX\#\sY;\sJ) \subset \sE/\sK(\sX\#\sY;\sJ)$, 
$\sE/\sK(V,\sX\#\sY;\sJ) \subset \sE/\sK(\sX\#\sY;\sJ)$ 
will correspond to in terms of $\sE/\sK(\sX\backslash \{x\},\sX;\sJ)$, $\sE/\sK(\sY\backslash\{y\},\sY;\sJ)$ and $\sE/\sK((0,1) \times S^{n-1},[0,1] \times S^{n-1};\sJ)$. In the end we get an exact sequence which we record as the next theorem.

\medskip
\noindent
{\bf Theorem 5.2.} {\em The Mayer--Vietoris exact sequence gives rise to an exact sequence
\[
\scriptsize{\begin{array}{c}
K_0(\sE/\sK((0,1) \times S^{n-1},[0,1] \times S^{n-1};\sJ)) \to K_0(\sE/\sK(\sX\backslash\{x\},\sX;\sJ)) \oplus K_0(\sE/\sK(\sY\backslash\{y\},\sY;\sJ) \to K_0(\sE/\sK(\sX\#\sY;\sJ)) \\
\uparrow\hskip 3.25 in \downarrow \\
K_1(\sE/\sK(\sX\#\sY;\sJ)) \leftarrow K_1(\sE/\sK(\sX\backslash \{x\},\sX_j;\sJ)) \oplus K_1(\sE/\sK(\sY\backslash\{y\},\sY;\sJ)) \leftarrow K_1(\sE/\sK((0,1) \times S^{n-1},[0,1] \times S^{n-1};\sJ))
\end{array}}
\]
}

\medskip
We begin the concluding remarks for this section with two $K$-theory (\cite{2},\cite{12}) consequences of Theorem~$5.1$.

If the dimension $n$ of $\sX$ is $\ge 3$, the homomorphism $\Psi: \sE/\sK(\sX;\sC^-_n) \to \pi(M(L^{\infty}(\sX,\mu)))$ so that $\Psi \mid \pi(M(L^{\infty}(\sX,\mu))) = id_{\pi(M(L^{\infty}(\sX,\mu)))}$ has immediately the following $K$-theory consequences.

\medskip
\noindent
{\bf Corollary 5.1.} {\em Assume $n = \dim \sX \ge 3$. Then $K_0(L^{\infty}(\sX;\mu)) \simeq \{f \in L^{\infty}(\sX,\mu) \mid f(\sX) \subset {\mathbb Z}\}$ is isomorphic to a direct summand of $K_0(\sE/\sK(\sX,\sJ))$ and the isomorphism is also w.r.t.\ the order structure. We also have that
\[
K_1(\sE/\sK(\sX;\sJ)) \simeq K_1(\ker \Psi).
\]
}

\medskip
By Thm.~$5.1$ if $\sJ \supseteq \sC^-_n$ we have that $\pi(M(C(\sX))) \simeq C(\sX)$ is the center of $\sE/\sK(\sX;\sJ)$. This yields a homomorphism $Z$
\[
C(\sX) \otimes \sE/\sK(\sX;\sJ) \to \sE/\sK(\sX;\sJ)
\]
where the LHS has a norm equivalent to that of the tensor product of the $C^*$-algebras $C(\sX)$ and $p(\sE(\sX;\sJ)))$ and which maps $f \otimes 1$ to $\pi(M(f))$ and $1 \otimes T$ to $T$. Moreover if $\Delta: C(\sX) \otimes C(\sX) \to C(\sX)$ is the restriction to the diagonal then $Z \circ (\Delta \otimes id) = Z \circ (id \otimes Z)$. This then can be used to get that $K_0(\sE/\sK(\sX;\sJ))$ is a $K_0(C(\sX))$-module and a similar reasoning with $C(\sX) \otimes C(S^1)$ and $\sE/\sK(\sX;\sJ) \otimes C(S^1)$ gives that the ${\mathbb Z}/2{\mathbb Z}$ graded $K_*(\sE/\sK(\sX;\sJ))$ is a $K_*(C(\sX))$-module, which is certainly not a new thing about the center of a $C^*$-algebra $p(\sE(\sX;\sC^-_n))$ so that we deal with $C^*$-algebras.

\medskip
\noindent
{\bf Corollary 5.2.} {\em $K_*(\sE/\sK(\sX;\sJ))$ is a $K_*(C(\sX))$-module when $\sJ \supseteq \sC^-_n$.
}

\medskip
\noindent
{\bf Remark 5.1.} The fact that Thm.~$5.1.$c) identifies the center of $\sE/\sK(\sX;\sJ)$ when $\sJ \supseteq \sC^-_n$ with $C(\sX)$ means that $\sE/\sK(\sX;\sJ)$ is a $C(\sX)$-algebra in the sense of G.~G. Kasparov \cite{16}, a notion which was generalized by E.~Kirchberg \cite{18} to so called $\sX$-algebra or algebra over a topological space. This suggests that when studying the $K$-theory of $\sE/\sK(\sX;\sJ)$ it may be natural to consider $K$-groups which take the $\sX$-structure into account (\cite{1}, \cite{18}, \cite{20}). Note however that with $\sE/\sK(\sX;\sJ)$ one will have to deal with a rather different class of $\sX$-algebras than those in (\cite{1}, \cite{7}, \cite{18}, \cite{20}), algebras which are non-separable and far from bundles. However, at least when $\sX$ has no boundary, the action of the diffeomorphism group of $\sX$ on $\sE/\sK(\sX;\sJ)$ provides some homogeneity.

\medskip
\noindent
{\bf Remark 5.2.} From these considerations about $\sE/\sK(\sX;\sJ)$ it seems that certain aspects with operations, like performing connected sums may be somewhat manageable. On the other hand in the simple cases when $\sX$ is a sphere or a ball or a product of spheres and balls there is a lot of mystery about the $\sE/\sK(\sX;\sJ)$. On the analysis side the smooth functions on $\sX$ which arise from the bicommutant $\mod \sJ$ construction in \S 4 is certainly an important analytic question about $\sE/\sK(\sX;\sJ)$.

\medskip
\noindent
{\bf Remark 5.3.} Since the smooth algebra $\sD(\tau,\sJ)$ depends only on $\sE(\tau;\sJ)$ we see that $\sD(\tau;\sJ)$ where $\tau$ is a $m$-tuple of multiplication operators from an embedding of $\sX$ into ${\mathbb R}^m$, does not depend on the choice of the embedding. Thus we get an algebra $\sD(\sX;\sJ)$ which depends only on $\sX$. If the dimension of $\sX$ is $> 0$ we also have that $\sD(\sX;\sJ) \cap \sK = \{0\}$. Thus $\sD(\sX;\sJ)$ and $\sD/\sK(\sX;\sJ)$ are isomorphic if $\dim \sX > 0$.



\begin{thebibliography}{99}

\bibitem{1} R.~Bentnam, {\em Kirchberg $X$-algebra with real rank zero and intermediate cancellation}, J. Noncommut. Geom. {\bf 8} (2014), 1061--1081.

\bibitem{2} B.~Blackadar, {\em $K$-theory for operator algebras}, 2nd Ed., Math. Sci. Res. Inst. Publ. 5, Cambridge Univ. Press, Cambridge, 1998. 

\bibitem{3} A.~Ber, J.~Huang, G.~Levitina, and F.~Sukochev, {\em Derivations with values in ideals of semifinite von~Neumann algebras}, J. Funct. Anal. {\bf 272} (2017), 4984--4997.

\bibitem{4} H.~Bercovici and D.~V. Voiculescu, {\em The analogue of Kuroda's theorem for $n$-tuples}, The Gohberg Anniversary Collection, Vol.~II (Calgary AB (1988), 57--60), Oper. Theory Adv. Appl. {\bf 41}, Birkhauser, Basel (1989).

\bibitem{5} J.~Bourgain and D.~V. Voiculescu, {\em The essential centre of the {\rm mod} a diagonalization ideal commutant of an $n$-tuple of commuting Hermitian operators}, Oper. Theory Adv. Appl. {\bf 252}, 77--80, Birkhauser/Springer (2016).

\bibitem{6} A.~Connes, {\em On the spectral characterization of manifolds}, J.~Noncommut. Geom {\bf 7} (2013), No.~1, 1--82.

\bibitem{7} M.~Dadarlat and W.~Winter, {\em Trivialization of $C(X)$-algebras with strongly self-absorbing fibres}, Bull. Soc. Math. France {\bf 136} (2008), No.~4, 575--606.

\bibitem{8} G.~David and D.~V. Voiculescu, {\em $s$-numbers of singular integrals for the invariance of absolutely continuous spectra in fractional dimension}, J. Funct. Anal. {\bf 94} (1990), 14--26.

\bibitem{9} W.~Donoghue, {\em On the perturbation of spectra}, Comm. Pure Appl. Math., Vol.~18 (1965), 559--579.

\bibitem{10} K.~J. Dykema and A. Skripka, {\em Perturbation formulas for traces on normed ideals}, Comm. Math. Phys., No.~3 (325):1107--1138 (2014).

\bibitem{11} I.~C. Gohberg and M.~G. Krein, {\em Introduction to the theory of non-selfadjoint operators}, Translations of Mathematical Monographs, Vol.~18, AMS, Providence, RI (2005).

\bibitem{12} N.~Higson and J.~Roe, ``Analytic $K$-homology,'' Oxford Mathematical Monographs, Oxford University Press, Oxford, 2000.

\bibitem{13} M.~Hirsch, ``Differential Topology,'' Springer Verlag, 1976.

\bibitem{14} B.~Johnson and S.~Parrott, {\em Operators commuting with a von~Neumann algebra modulo the set of compact operators}, J. Funct. Anal. {\bf 11} (1972), 39--61.

\bibitem{15} V.~Kaftal and G.~Weiss, {\em Compact derivations relative to semifinite von~Neumann algebras}, J. Funct. Anal. {\bf 62} (1985), 202--220.

\bibitem{16} G.~G. Kasparov, {\em Equivariant $KK$-theory and the Novikov conjecture}, Invent. Math. {\bf 91} (1988), 147--201.

\bibitem{17} E.~Kissin, D.~Potapov, V.~Shulman, and F.~Sukochev, {\em Operator smoothness in Schatten norms for functions of several variables: Lipschitz conditions, differentiability and unbounded derivations}, Proc. Lond. Math. Soc. (3), No.~4 (105), (2012), 661--702.

\bibitem{18} E.~Kirchberg, {\em Das nicht-kommutative Michael--Auswahlprinzip und die Klassifikation nicht-zinfacher Algebren, $C^*$-algebras (M\"{u}nster 1999)}, 92--141, Springer, Berlin, 2000.

\bibitem{19} C.~Liaw and S.~Treil, {\em Singular integrals rank one perturbations and Clark model in general situation, Harmonic analysis, partial differential equations, Banach spaces and operator theory}, Vol.~2, 85--132, Assoc. Women Math., Ser.~5, Springer Cham., 2017.

\bibitem{20} R.~Meyer and R.~Nest, {\em $C^*$-Algebras over Topological Spaces: Filtrated $K$-theory}, Canad. J. Math., Vol.~64(2), pp.~368--408 (2012).

\bibitem{21} F.~L. Nazarov and V.~V. Peller, {\em Functions of perturbed $n$-tuples of commuting self-adjoint operators}, J. Funct. Anal. {\bf 266} (2014), 5398--5428.

\bibitem{22} A.~G. Poltoratski, {\em Equivalence up to a rank one perturbation}, Pacific J. Math. {\bf 194} (2000), No.~1, 175--188.

\bibitem{23} S.~Popa, {\em The commutant modulo the set of compact operators of a von~Neumann algebra}, J. Funct. Anal. {\bf 71} (1987), 393--408.

\bibitem{24} B.~Simon, {\em Trace ideals and their applications}, 2nd Ed., Mathematical Surveys and Monographs {\bf 120}, AMS; Providence, RI, 2005.

\bibitem{25} B.~Simon and T.~Wolff, {\em Singular continuous spectrum under rank one perturbations and localization for random Hamiltonians}, Comm. Pure Appl. Math. {\bf 39} (1986), No.~1, 75--90.

\bibitem{26} D.~V. Voiculescu, {\em Some results on norm-ideal perturbations of Hilbert space operators I}, J. Operators Theory {\bf 2} (1979), 3--37.

\bibitem{27} D.~V. Voiculescu, {\em Some results on norm-ideal perturbations of Hilbert space operators II}, J. Operator Theory {\bf 5} (1981), 77--100.

\bibitem{28} D.~V. Voiculescu, {\em On the existence of quasicentral approximate units relative to normed ideals I}, J. Funct. Anal. {\bf 91} (1)(1991), 1--36.

\bibitem{29} D.~V. Voiculescu, {\em $K$-theory and perturbations of absolutely continuous spectra}, Comm. Math. Phys. {\bf 365} (2019), No.~1, 363--373.

\bibitem{30} D.~V. Voiculescu, {\em Some results and a $K$-theory problem about threshold commutants {\rm mod} normed ideals}, (preprint) arXiv:1911.07377.

\bibitem{31} D.~V. Voiculescu, {\em Commutants {\rm mod} normed ideals}, in ``Advances in Non-commutative Geometry On the Occasion of Alain Connes' 70th Birthday'', Springer Verlag (2020), 585--606.

\bibitem{32} D.~V. Voiculescu, {\em The formula for the quasicentral modulus in the case of spectral measures on fractals}, (preprint) arXiv:2006.14456.

\end{thebibliography}
\end{document}